\documentclass{elsart}
\usepackage{ae}
\usepackage{graphicx}
\usepackage{amsmath}
\usepackage{amssymb}
\usepackage{amsfonts}
\usepackage{latexsym}
\usepackage{verbatim}
\usepackage{url}
\newcommand{\cB}{{\mathcal B}}

\newcommand{\cA}{{\mathcal A}}

\newcommand{\cF}{{\mathcal F}}

\newcommand{\NN}{{\mathbb N}}

\newcommand{\cAstar}{\mathcal{A}^*}
\newcommand{\cAplus}{\mathcal{A}^+}
\newcommand{\empt}{\varepsilon}
\newcommand{\cAw}{\mathcal{A}^{\omega}}
\newcommand{\cAinf}{\mathcal{A}^{\infty}}

\newcommand{\rev}{\widetilde}

\newcommand{\bx}{\mathbf{x}}

\newcommand{\pref}{\prec_p}
\newcommand{\suff}{\prec_s}
\newcommand{\bs}{\mathbf{s}}
\newcommand{\bt}{\mathbf{t}}

\newcommand{\by}{\mathbf{y}}
\newcommand{\bv}{\mathbf{v}}
\newcommand{\bu}{\mathbf{u}}
\newcommand{\bz}{\mathbf{z}}

\newcommand{\bbf}{\mathbf{f}}

\newcommand{\Ult}{\mbox{Ult}}

\numberwithin{thm}{section}

\newenvironment{remq}{\begin{rem}\rm }{\end{rem}}
\newenvironment{example}{\begin{exmp}\rm }{\end{exmp}}
\numberwithin{thm}{section}

\journal{Theoretical Computer Science}

\begin{document}

% -----------------------Frontmatter-----------------------------

\begin{frontmatter}

\title{A characterization of fine words over a finite alphabet}

%\thanks[CIL]{Research supported by CRM, ISM, and LaCIM.}

\author{Amy Glen%\thanksref{CIL}
}

\address{LaCIM, Universit\'e du Qu\'ebec \`a Montr\'eal, C.P. 8888, succursale Centre-ville, Montr\'eal, Qu\'ebec, CANADA, H3C 3P8}
\ead{amy.glen@gmail.com}

\date{March 20, 2007}

\begin{abstract}
To any infinite word $\bt$ over a finite alphabet $\cA$ we can associate two infinite words $\min(\bt)$ and $\max(\bt)$ such that any prefix of $\min(\bt)$ (resp.~$\max(\bt)$) is the \emph{lexicographically} smallest (resp.~greatest)  amongst the factors of $\bt$ of the same length. We say that an infinite word $\bt$ over $\cA$ is {\em fine} if  there exists an infinite word $\bs$ such that, for any lexicographic order, $\min(\bt) = a\bs$ where $a = \min(\cA)$. In this paper, we characterize fine words; specifically, we prove that an infinite word $\bt$ is fine if and only if $\bt$ is either a \emph{strict episturmian word} or a strict ``skew episturmian word''. This characterization generalizes a recent result of G. Pirillo, who proved that a fine word over a 2-letter alphabet is either an (aperiodic) \emph{Sturmian word}, or an ultimately periodic (but not periodic) infinite word, all of whose factors are (finite) Sturmian.
\end{abstract}

\end{frontmatter}

\vspace*{-1cm}\section{Introduction}

To any infinite word $\bt$ over a finite alphabet $\cA$ we can associate two infinite words $\min(\bt)$ and $\max(\bt)$ such that any prefix of $\min(\bt)$ (resp.~$\max(\bt)$) is the \emph{lexicographically} smallest (resp.~greatest)  amongst the factors of $\bt$ of the same length (see Pirillo \cite{gP05ineq}). In the recent paper \cite{gP05mors}, Pirillo defined \emph{fine words} over two letters; specifically, an infinite word $\bt$ over a 2-letter alphabet $\{a,b\}$ ($a < b$) is said to be \emph{fine} if  $(\min(\bt), \max(\bt)) = (a\bs, b\bs)$ for some infinite word $\bs$. Pirillo \cite{gP05mors} characterized these words, and remarked that perhaps his characterization can be generalized to an arbitrary finite alphabet; we do just that in this paper. Firstly, we extend the definition of a fine word to more than two letters. That is, we say that an infinite word $\bt$ over $\cA$ is \emph{fine} if  there exists an infinite word $\bs$ such that, for any lexicographic order, $\min(\bt) = a\bs$  where $a = \min(\cA)$. Roughly speaking, our main result states that an infinite word $\bt$ is fine if and only if $\bt$ is either a \emph{strict episturmian word} or a strict ``skew episturmian word'' (i.e., a particular kind of non-recurrent infinite word, all of whose factors are {\em finite episturmian}).

\section{Notation and terminology}

{\bf Finite and infinite words}

Let $\cA$ denote a finite alphabet. A (finite) \emph{word} over $\cA$ is an
element of the \emph{free monoid} $\cAstar$ generated by $\cA$, in
the sense of concatenation. The identity $\empt$ of $\cAstar$ is
called the \emph{empty word}, and the \emph{free semigroup},
denoted by $\cAplus$, is defined by $\cAplus :=
\cAstar\setminus\{\empt\}$.

Given $w = x_{1}x_{2}\cdots x_{m} \in \cAplus$ with each 
$x_{i} \in \cA$, the \emph{length} of $w$ is $|w| = m$  
 (note that $|\empt| = 0$). The \emph{reversal} $\rev{w}$ of $w$  
is given by $\rev{w} = x_{m}x_{m-1}\cdots x_{1}$, and if 
$w = \rev{w}$, then $w$ is called a \emph{palindrome}.

An \emph{infinite word} (or simply \emph{sequence}) $\bx$ is a sequence indexed by $\NN$ with values in $\cA$, i.e., $\bx = x_0x_1x_2\cdots$ with each $x_i \in \cA$. The set of 
all infinite words over $\cA$ is denoted by $\cAw$, and we define $\cAinf := \cAstar \cup \cAw$. 
An {\em ultimately periodic} infinite word can be written as $uv^\omega = uvvv\cdots$, for some $u$, $v \in \cAstar$, $v\ne \empt$. If $u = \empt$, then such a word is {\em periodic}. 
An infinite word that is not ultimately periodic is said to be {\em aperiodic}.

A finite word $w$ is a \emph{factor} of $z \in \cAinf$ if $z =
uwv$ for some $u \in \cAstar$, $v \in \cAinf$.  
Further, $w$ is called a \emph{prefix} (resp.
\emph{suffix}) of $z$ if $u = \empt$ (resp.~$v = \empt$), and we
write $w \pref z$ (resp.~$w \suff z$).  We say that $w$ is a \emph{proper} factor (resp.~prefix, suffix) of $z$ if $uv \ne \empt$ (resp.~$v\ne\empt$, $u\ne \empt$).  
An infinite word $\bx \in \cAw$ is called a \emph{suffix} of $\bz \in \cAw$ if there exists a word $w \in \cA^*$ such that $\bz = w\bx$.   
A factor $w$ of a
word $z \in \cAinf$ is \emph{right} (resp.~\emph{left})
\emph{special} if $wa$, $wb$ (resp.~$aw$, $bw$) are factors of $z$
for some letters $a$, $b \in \cA$, $a \ne b$.

For $\bx \in \cAw$, $F(\bx)$ denotes the set of all its
factors, and $F_n(\bx)$ denotes the set of all factors of
$\bx$ of length $n \in \NN$, i.e., $F_n(\bx) := F(\bx)
\cap \cA^n$. Moreover, the \emph{alphabet} of $\bx$ is Alph$(\bx)
:= F(\bx) \cap \cA$, and we denote by Ult$(\bx)$ the set of
all letters occurring infinitely often in $\bx$.  Any two infinite words $\bx$, $\by \in \cAw$ are said to be \emph{equivalent} if $F(\bx) = F(\by)$, i.e., if $\bx$ and $\by$ have the same set of factors. A factor of an infinite word $\bx$ is \emph{recurrent} in $\bx$ if it occurs infinitely many times in $\bx$, and $\bx$ itself is said to be \emph{recurrent} if all of its factors are recurrent in it.

{\bf Lexicographic order}

Suppose the alphabet $\cA$ is totally ordered by the relation $<$. Then we 
can totally order $\cA^*$ by the \emph{lexicographic order} $<$, 
defined as follows. Given two words $u$, $v \in \cA^+$, we have $u
< v$ if and only if either $u$ is a proper prefix of $v$ or $u =
xau^\prime$ and $v = xbv^\prime$, for some $x$, $u^\prime$,
$v^\prime \in \cAstar$ and letters $a$, $b$ with $a < b$. This is
the usual alphabetic ordering in a dictionary, and we say that $u$
is \emph{lexicographically less} than $v$. This notion 
naturally extends to $\cAw$, as follows. Let $\bu =
u_0u_1u_2\cdots$ and $\bv = v_0v_1v_2\cdots$, where $u_j$, $v_j
\in \cA$. We define $\bu < \bv$ if there exists an index $i\geq0$
such that $u_j = v_j$ for all $j=0,\ldots, i-1$ and $u_{i} < v_{i}$. Naturally, $\leq$ will mean $<$ or $=$.

Let $w \in \cAinf$ and let $k$ be a positive integer. We denote by $\min(w | k)$ (resp.~$\max(w | k)$) the lexicographically smallest (resp.~greatest) factor of $w$ of length $k$ for the given order (where $|w|\geq k$ for $w$ finite). If $w$ is infinite, then it is clear that $\min(w | k)$ and $\max(w | k)$ are prefixes of the respective words $\min(w | k+1)$ and $\max(w | k+1)$. So we can define, by taking limits, the following two infinite words (see \cite{gP05ineq})
\[
  \min(w) = \underset{k\rightarrow\infty} {\lim}\min(w | k) \quad \mbox{and} \quad 
  \max(w) = \underset{k\rightarrow\infty}{\lim}\max(w | k).
\]

{\bf Morphisms and the free group}

A \emph{morphism on} $\cA$ is a map $\psi: \cAstar \rightarrow
\cAstar$ such that $\psi(uv) = \psi(u)\psi(v)$ for all $u, v \in
\cAstar$.  It is uniquely determined by its image on the alphabet
$\cA$. All morphisms considered in this paper will be non-erasing: the image of any non-empty word is never empty. Hence the action of a morphism $\psi$ on $\cA^*$ naturally extends to infinite words; that is, if $\bx = x_0x_1x_2 \cdots \in \cAw$, then $\psi(\bx) = \psi(x_0)\psi(x_1)\psi(x_2)\cdots$. 

The free monoid $\cAstar$ can be naturally embedded within a free group. We denote by 
$\cF(\cA)$ the free group over $\cA$ that properly contains $\cA$, and is 
obtained from $\cA$ by adjoining the inverse $a^{-1}$ of each letter 
$a \in \cA$. More precisely, we construct a new alphabet $\cA^{\pm}$ that 
consists of all letters $a$ of $\cA$ and their `inverses' $a^{-1}$, i.e., 
$\cA^{\pm} = \{a, ~a^{-1} ~|~ a \in \cA\}$.  If one defines on the free 
monoid $(\cA^{\pm})^*$ the involution $(a^{-1})^{-1} = a$ 
for each $a \in (\cA^{\pm})^*$, then necessarily, we 
have $(uv)^{-1} = v^{-1}u^{-1}$ for all $u$, $v \in (\cA^\pm)^*$. The 
free group $\cF(\cA)$ over $\cA$ is the quotient of $(\cA^\pm)^*$ under 
the relation: $aa^{-1} = a^{-1}a = \empt$ for all $a \in \cA$. In what follows, we use the notation $p^{-1}w$ and $ws^{-1}$ to indicate the removal of a prefix $p$ (resp.~suffix $s$) from a finite word $w$.

Any morphism $\psi$ on $\cA$ can be uniquely extended to an endomorphism of 
$\cF(\cA)$ by defining $\psi(a^{-1}) = (\psi(a))^{-1}$ for each 
$a \in \cA$. %, from which it follows that $\sigma(w^{-1}) = (\sigma(w))^{-1}$ for any $w \in \cF(\cA)$.

\section{Episturmian words}

An interesting generalization of \emph{Sturmian words} (i.e., aperiodic infinite words of minimal complexity) to a finite alphabet is the family of \emph{Arnoux-Rauzy sequences}, the study of which began in  \cite{pAgR91repr} (also see \cite{ jJgP02onac, rRlZ00agen} for example). More recently, a slightly wider class of infinite words, aptly called \emph{episturmian words}, was introduced by Droubay, Justin, and Pirillo \cite{xDjJgP01epis} (also see  \cite{aG05powe, jJgP02epis, jJgP04epis, jJlV00retu} for instance). 
An infinite word $\bt \in \cAw$ is \emph{episturmian} if $F(\bt)$ is closed under
reversal and $\bt$ has at most one right (or equivalently left) special factor of each length. Moreover, an episturmian word is \emph{standard} if all of its left special factors are prefixes of
it. Sturmian words are exactly the aperiodic episturmian words over a 2-letter alphabet. 

Standard episturmian words were characterized in \cite{xDjJgP01epis} using the concept of the \emph{palindromic right-closure} $w^{(+)}$ of a finite word $w$, which is the (unique) shortest palindrome having $w$ as a prefix (see \cite{aD97stur}).   
Specifically, an infinite word $\bs \in \cAw$ is standard episturmian if and only if  there exists an infinite word $\Delta(\bs) = x_1x_2x_3\ldots$ ($x_i
\in \cA$), called the \emph{directive word} of $\bs$, such that the infinite sequence of palindromic prefixes $u_1 =
\empt$, $u_2$, $u_3$, $\ldots$ of $\bs$ (which exists by results in \cite{xDjJgP01epis}) is given by 
\begin{equation} \label{eq:02.09.04}
  u_{n+1} = (u_nx_n)^{(+)}, \quad 
  n \in \NN^+.
\end{equation}
This characterization extends to the case of an arbitrary finite alphabet a construction given in \cite{aD97stur} for all {\em standard Sturmian words}. An important point is that a standard episturmian word 
$\bs$ can be constructed as a limit of an infinite sequence of its palindromic prefixes, i.e., 
$\bs = \lim_{n\rightarrow\infty} u_n$. 

\noindent{\it Note.} Episturmian words are (uniformly) recurrent \cite{xDjJgP01epis}.

\subsection{Relation with episturmian morphisms}

Let $a \in \cA$ and denote by $\Psi_a$ the morphism on $\cA$ defined by 
\[
\Psi_a : \left\{\begin{array}{lll}
                a &\mapsto &a \\
                x &\mapsto &ax \quad \mbox{for all $x \in \cA\setminus\{a\}$}.
                \end{array} \right.
\]
Together with the permutations of the alphabet, all of the morphisms $\Psi_a$ generate by composition the monoid of {\em epistandard morphisms} (`epistandard' is an elegant shortcut for `standard episturmian' due to Richomme \cite{gR03conj}). The submonoid generated by the $\Psi_a$ only is the monoid of {\em pure epistandard morphisms}, which includes the {\em identity morphism} Id$_{\cA} =$ Id, and consists of all the \emph{pure standard (Sturmian)  morphisms} when $|\cA|=2$. 

When viewed as an endomorphism of the free group $\cF(\cA)$, the morphism $\Psi_a$ is invertible; that is, $\Psi_a$ is a {\em positive automorphism} of $\cF(\cA)$, and its inverse is given by
\[
\Psi_a^{-1} : \left\{\begin{array}{lll}
                a &\mapsto &a \\
                x &\mapsto &a^{-1}x \quad \mbox{for all $x \in \cA\setminus\{a\}$}.
                \end{array} \right.
\]
It follows that every epistandard morphism is a (positive) automorphism of $\cF(\cA)$. See \cite{aG06onst, eG06auto, gR03conj, zWyZ99some} for work involving the invertibility of episturmian morphisms.   

\begin{remq} If $\bx = \Psi_a(\by)$ or $\bx = a^{-1}\Psi_a(\by)$ for some $\by \in \cAw$ and $a \in \cA$, then the letter $a$ is \emph{separating for $\bx$} and its factors; that is, any factor of $\bx$ of length 2 contains the letter $a$.
\end{remq}

Another useful characterization of standard episturmian words is the following (see \cite{jJgP02epis}). An infinite word $\bs \in \cAw$ is standard episturmian with directive word $\Delta(\bs) = x_1x_2x_3\cdots$ ($x_i \in \cA$) if and only if there exists an infinite sequence of recurrent infinite words
$\bs^{(0)} = \bs$, $\bs^{(1)}$, $\bs^{(2)}$, $\ldots$ such that $\bs^{(i-1)} = \Psi_{x_i}(\bs^{(i)})$ for all $i \in \NN^+$. Moreover, each $\bs^{(i)}$ is a standard episturmian word with directive word 
$\Delta(\bs^{(i)}) = x_{i+1}x_{i+2}x_{i+3}\cdots$, the \emph{$i$-th shift} of $\Delta(\bs)$. 

To the prefixes of the directive word $\Delta(\bs) = x_1x_2\cdots$, we associate the morphisms 
\[
  \mu_0 := \mbox{Id}, \quad \mu_n := \Psi_{x_1}\Psi_{x_2}\cdots\Psi_{x_n}, \quad n \in \NN^+, 
\]
and define the words 
\[
  h_n := \mu_n(x_{n+1}), \quad n \in \NN, 
\]
which are clearly prefixes of $\bs$.  For the palindromic prefixes $(u_i)_{i\geq1}$ given by \eqref{eq:02.09.04}, we have the following useful formula \cite{jJgP02epis}
\[
  u_{n+1} = h_{n-1}u_{n};
\]
whence, for $n > 1$ and $0 < p < n$, 
\begin{equation} \label{eq:u_n}
  u_n = h_{n-2}h_{n-3}\cdots h_1h_0 = h_{n-2}h_{n-3}\cdots h_{p-1}u_p.
\end{equation}

\begin{remq} 
Evidently, if a standard episturmian word $\bs$ begins with the letter $x \in \cA$, then $x$ is separating for $\bs$ (see \cite[Lemma 4]{xDjJgP01epis}).
\end{remq}

\subsection{Strict episturmian words}

A standard episturmian word $\bs\in \cA^\omega$, or any equivalent (episturmian) word,   
is said to be \emph{$\cB$-strict} (or $k$-\emph{strict} if $|\cB|=k$, or {\em strict} if $\cB$ is understood) if 
Alph$(\Delta(\bs)) =$ Ult$(\Delta(\bs)) = \cB \subseteq \cA$.  In particular, a standard episturmian word 
over $\cA$ is $\cA$-strict if every letter in $\cA$ occurs infinitely often in its directive word. 
The $k$-strict episturmian words have complexity $(k - 1)n + 1$ for 
each $n \in \NN$; such words are exactly the $k$-letter Arnoux-Rauzy  
sequences. %The strict standard episturmian words are precisely the standard (or characteristic) Arnoux-Rauzy sequences.  
Note that the $2$-strict episturmian words correspond to the (aperiodic) Sturmian words.

\begin{remq}\label{R:ult_remark}
Suppose $\bs \in \cA^\omega$ is a standard episturmian word. If $\bs$ is not $\cA$-strict, then $\Ult(\Delta(\bs)) = \cB \subset \cA$ and there exists a $\cB$-strict standard episturmian word $\bs'$ and a pure epistandard morphism $\mu$ on $\cA$  such that $\bs = \mu(\bs')$. More precisely, let $\Delta(\bs) = x_1x_2x_3\cdots$ and let $m$ be minimal such that Alph$(x_{m+1}x_{m+2}\cdots) = \cB \subset \cA$. That is, $x_1x_2\cdots x_m$ is the shortest prefix of $\Delta(\bs)$ that contains all the letters not appearing infinitely often in $\Delta(\bs)$, namely the letters in $\cA \setminus \cB$. Then 
$\bs = \mu_m(\bs^{(m)})$ where $\bs^{(m)}$ is the $\cB$-strict standard episturmian word with directive word $\Delta(\bs^{(m)}) = x_{m+1}x_{m+2}\cdots$. For example, if $\Delta(\bs) = c(ab)^\omega$, then $\bs = \Psi_c(\bs^{(1)})$ where $\Delta(\bs^{(1)}) = (ab)^\omega$, i.e., $\bs^{(1)}$ is the well-known \emph{Fibonacci word} over $\{a,b\}$. 
\end{remq}

\section{Fine words}

Recall that an infinite word $\bt$ over $\cA$ is \emph{fine} if  there exists an infinite word $\bs$ such that, for any lexicographic order, $\min(\bt) = a\bs$ where $a = \min(\cA)$.

\noindent{\it Note.} Since there are only two lexicographic orders on words over a 2-letter alphabet, a fine word $\bt$ over $\{a,b\}$ ($a< b$) satisfies $(\min(\bt), \max(\bt)) = (a\bs, b\bs)$ for some infinite word $\bs$. 

Recently, Pirillo \cite{gP05mors} characterized fine words over a 2-letter alphabet. Specifically:

\begin{prop}  \label{P:gP05mors} {\em \cite{gP05mors}} Suppose $\bt$ is an infinite word over $\{a,b\}$. Then the following properties are equivalent:  
\begin{enumerate}
\item[i)] $\bt$ is fine;
\item[ii)] either $\bt$ is a Sturmian word, or $\bt = v\mu(x)^\omega$ where $\mu$ is a pure standard Sturmian morphism on $\{a,b\}$, and $v$ is a non-empty suffix of $\mu(x^py)$ for some $p \in \NN$ and $x$, $y \in \{a,b\}$ $(x\ne y)$.  
\end{enumerate}   
\end{prop}

In other words, a fine word over two letters is either a Sturmian word or an ultimately periodic (but not periodic) infinite word, all of whose factors are (finite) Sturmian, i.e., a so-called \emph{skew Sturmian} word (see \cite{gHmM40symb}). In this paper, we generalize Pirillo's result to infinite words over two or more  letters. 

The next two propositions are needed for the proof of our main result (Theorem \ref{T:fine}, to follow). Recall that the Arnoux-Rauzy sequences are precisely the strict episturmian words.

\begin{prop} \label{P:jJgP-1} \emph{\cite{jJgP02onac}} Suppose $\bs$ is an infinite word over a finite alphabet $\cA$. Then the following properties are equivalent:  
\begin{enumerate}
\item[i)] $\bs$ is a standard Arnoux-Rauzy sequence;
\item[ii)] $a\bs = \min(\bs)$ for any letter $a \in \cA$ and lexicographic order $<$ satisfying $a = \min(\cA)$. 
\end{enumerate}   
\end{prop}

\begin{prop}  \label{P:gP-2} \emph{\cite{gP05ineq}} Suppose $\bs$ is an infinite word over a finite alphabet $\cA$. Then the following properties are equivalent:   
\begin{enumerate}
\item[i)] $\bs$ is standard episturmian;
\item[ii)] $a\bs \leq \min(\bs)$ for any letter $a \in \cA$ and lexicographic order $<$ satisfying $a = \min(\cA)$. 
\end{enumerate} 
\end{prop}

The following key lemma is also needed. From now on, it will be convenient to denote by $\bv_p$ the prefix of length $p$ of a given infinite word $\bv$.

\begin{lem} \label{L:25.01.2006} Let $\cA$ be a finite alphabet and let $a \in \cA$.  
Suppose $\bt$, $\bs \in \cAw$ are infinite words such that 
$\bt = \Psi_z(\bt^{(1)})$ and $\bs = \Psi_z(\bs^{(1)})$ for some $z \in \mbox{{\em Alph}}(\bt^{(1)})$.  Then 
\[
  \min(\bt^{(1)})  = a\bs^{(1)}  \quad \Leftrightarrow \quad \min(\bt) = a\bs. 
\]
\end{lem}

\begin{remq} \label{R:more_general} Let $\bt$, $\bt^{(1)}$, $\bs$, $\bs^{(1)} \in \cA^\omega$ be such that 
$\bt = \Psi_z(\bt^{(1)})$ and $\bs = \Psi_z(\bs^{(1)})$ for some letter $z$ (not necessarily in $\mbox{Alph}(\bt^{(1)})$).  Using  similar reasoning as in the proof below, it can be shown that 
\[
 \min(\bt^{(1)})  = a\bs^{(1)} \quad \Leftrightarrow \quad 
 \min(\bt) = \begin{cases} 
  za\bs  &\mbox{if $z < a$}, \\
  a\bs &\mbox{if $z \geq a$}.
  \end{cases}
\]  
For example, let $\cA = \{a,b,c\}$ with $a<b<c$ and suppose $\bbf$ is the Fibonacci word over $\{a,b\}$  (i.e., the standard episturmian word directed by $(ab)^\omega$). Then $\min(\bbf) =  a\bbf$, and hence $\min(\Psi_c(\bbf)) = a\Psi_c(\bbf)$. On the other hand, if $\bbf^\prime$ is the Fibonacci word over $\{b,c\}$, then $\min(\bbf^\prime) = b\bbf^\prime$ and we have  $\min(\Psi_a(\bbf^\prime))  = ab\Psi_a(\bbf^\prime)$.  Lemma \ref{L:25.01.2006} is a special case of this result with $z \in \mbox{Alph}(\bt^{(1)}) \subseteq \cA$ and is sufficient for our purposes. 
\end{remq}

\begin{pf*}{Proof of Lemma $\ref{L:25.01.2006}$}

($\Leftarrow$):  We have $\min(\bt) = a\bs$. First observe that $a \in$ Alph$(\bt^{(1)})$. Indeed,  if $a = z$, then $a \in$ Alph$(\bt^{(1)})$ since $z \in \mbox{Alph}(\bt^{(1)})$ (in fact, $zz \in F(\bt)$ since $zz$ is a prefix of $a\bs = z\bs$, and hence $a = z \in$ Alph$(\bt^{(1)})$). On the other hand, if $a \ne z$, then we must have $a \in$ Alph$(\bt^{(1)})$, otherwise $a$ is not in the alphabet of $\bt = \Psi_z(\bt^{(1)})$, which is impossible since $F(a\bs) \subseteq F(\bt)$.

Now we show that 
$F(a\bs^{(1)}) \subseteq F(\bt^{(1)})$.  Suppose not, i.e., suppose 
$F(a\bs^{(1)}) \not\subseteq F(\bt^{(1)})$. Then there exists a minimal $m \in \NN^+$ such that $a\bs_{m}^{(1)} \not\in F(\bt^{(1)})$. Therefore, if $\bs_{m}^{(1)} = \bs_{m-1}^{(1)}x$ where $x\in \cA$, then 
\[
 a\bs_{m-1}^{(1)}x \not\in F(\bt^{(1)}). 
 \]
 Letting $\bs_{l} = \Psi_z(\bs_{m-1}^{(1)})$, we have 
 \[
a\Psi_z(\bs_{m-1}^{(1)}x) = a\bs_{l}\Psi_z(x) \in F(a\bs) \subseteq F(\bt),
\]
and hence $\Psi_z(a)\bs_{l}\Psi_z(x) \in F(\bt)$ since $z$ is separating for $\bt$. So, if $x \ne z$ then $a\bs_{m-1}^{(1)}x \in F(\bt^{(1)})$, which is impossible; whence $x = z$. But then,   
$\bs_{m+1}^{(1)} = \bs_{m-1}^{(1)}zy^\prime$ for some $y^\prime \in \cA$ and we have 
\[
 a \Psi_z(\bs_{m+1}^{(1)}) = a\bs_{l}z\Psi_z(y^\prime) \in F(a\bs) \subseteq F(\bt). 
 \]
 Thus, 
 \[
 \Psi_z(a)\bs_{l}zz \in F(\bt), 
 \]
and hence $a\bs_{m-1}^{(1)}z ~(= a\bs_{m-1}^{(1)}x)$ is a factor of $\bt^{(1)}$; a contradiction. Therefore, we conclude that $F(a\bs^{(1)}) \subseteq F(\bt^{(1)})$. 

Now suppose on the contrary that $\min(\bt^{(1)}) \ne  a\bs^{(1)}$. Then there exists a word $w^{(1)} \in F(\bt^{(1)})$ of minimal length $|w^{(1)}| = m$ such that 
\[
  w^{(1)} < a\bs_{m-1}^{(1)}. 
\]
Let $w^{(1)} = u^{(1)}x$ $(x \in \cA$) where $u^{(1)}$ is non-empty since $a \in$ Alph$(t^{(1)})$. Then, by minimality of $m$, $u^{(1)} \geq a\bs_{m-2}^{(1)}$, and therefore 
$u^{(1)} = a\bs_{m-2}^{(1)}$ (otherwise $u^{(1)} > a\bs_{m-2}^{(1)}$ implies $w^{(1)} > a\bs_{m-1}^{(1)}$). Hence,  
\[
  w^{(1)} = a\bs_{m-2}^{(1)}x \quad \mbox{with $w^{(1)} < a\bs_{m-1}^{(1)}$}, 
\]
and therefore
\[
  \bs_{m-1}^{(1)} = \bs_{m-2}^{(1)}y \quad \mbox{for some $y \in \cA$, $y>x$.}
\]
Now, letting $w = \Psi_z(w^{(1)})$ and $\bs_{l} = \Psi_z(\bs_{m-2}^{(1)})$, we have  
\[
  w = \Psi_z(a)\bs_{l}\Psi_z(x) \in F(\bt) \quad \mbox{and} \quad \bs_{l}\Psi_z(y) \pref \bs.
\]

Now consider $x^\prime$, $y^\prime \in \cA$ such that $w^{(1)}x^\prime \in F(\bt^{(1)})$ and 
$\bs_{m}^{(1)} = \bs_{m-2}^{(1)}yy^\prime$ is a prefix of $\bs^{(1)}$. Then $\Psi_z(w^{(1)}x^\prime)$ is a factor of $\bt$, where 
\[
\Psi_z(w^{(1)}x^\prime) = w\Psi_z(x^\prime) 
                                           =\Psi_z(a)\bs_l\Psi_z(x)\Psi_z(x^\prime) 
                                           = \begin{cases} 
                                                 \Psi_z(a)\bs_lz\Psi_z(x^\prime) &\mbox{if $x=z$}, \\
                                                 \Psi_z(a)\bs_lzx\Psi_z(x^\prime) &\mbox{if $x \ne z$}.
                                                 \end{cases}
\]
Therefore, the word $v = a\bs_lzx$ is a factor of $\bt$. Moreover,  $\Psi_z(\bs_{m}^{(1)})$ is a prefix of $\bs$, where 
\[
  \Psi_z(\bs_{m}^{(1)}) = \begin{cases} 
                                           \bs_lz\Psi_z(y^\prime) &\mbox{if $y = z$}, \\
                                           \bs_lzy\Psi_z(y^\prime) &\mbox{if $y \ne z$}, 
                                           \end{cases}
\]
and hence $\bs_{l+2} = \bs_lzy$ is a prefix of $\bs$. Accordingly, $v < a\bs_{l+2}$ since $x < y$, contradicting the fact that the prefixes of $a\bs$ are the lexicographically smallest factors of $\bt$. Thus,  we conclude that $\min(\bt^{(1)})  = a\bs^{(1)}$.

\noindent ($\Rightarrow$): We have $\min(\bt^{(1)}) = a\bs^{(1)}$. As above, it is easily shown that $F(a\bs) \subseteq F(\bt)$; whence $\min(\bt) \leq a\bs$.  Let us suppose $\min(\bt) \ne a\bs$. Then there exists a word $w \in F(\bt)$ of minimal length $|w| = l$ such that $w < a\bs_{l-1}$ 
If we let $w = ux$, $x \in \cA$, then $u\geq a\bs_{l-2}$, and hence $u = a\bs_{l-2}$ (otherwise $w > a\bs_{l-1}$). Therefore, 
\[
  w = a\bs_{l-2}x < a\bs_{l-1} 
\]
and hence $\bs_{l-1} = \bs_{l-2}y$ where $y \in \cA$, $y > x$.

Since the letter $z$ is separating for $\bt$, $\bs_{l-2}$ must end with $z$; otherwise $x = y = z$, which is impossible. Thus,  
\[
 w = a\bs_{l-3}zx \quad \mbox{and} \quad \bs_{l-1} = \bs_{l-3}zy, \quad y> x.
\]
Let $\bs_{m-1}^{(1)} = \Psi_z^{-1}(\bs_{l-1})$. We distinguish two cases: $y = z$ and $y\ne z$. 

\noindent\emph{Case $1$}: $y  = z$. We have $\bs_{l-1} = \bs_{l-3}zz$, and thus $\bs_{m-1}^{(1)}$ and $\bs_{m-2}^{(1)}$ both end with the letter $z$. Note that $z \ne a$ because $a \leq x < y = z$. Therefore, since $w$ begins with $a$ and $z$ is separating for $\bt$, we have $zw \in F(\bt)$.  
Now, 
\begin{align*}
  \Psi_z^{-1}(zw) &= \Psi_z^{-1}(za\bs_{l-1}z^{-1}x) \\
                             &=a\bs_{m-1}^{(1)}z^{-1}\Psi_z^{-1}(x) \\
                            &= a\bs_{m-1}^{(1)}z^{-1}z^{-1}x \\
                            &= a\bs_{m-3}^{(1)}x, 
\end{align*}
i.e., $zw = \Psi_z(w^{(1)})$ where $w^{(1)} = a\bs_{m-3}^{(1)}x \in F(\bt^{(1)})$. Therefore, as $\bs_{m-2}^{(1)}$ ends with $z > x$, we have $w^{(1)}  < a\bs_{m-2}^{(1)}$; a contradiction.  

\noindent\emph{Case $2$}: $y\ne z$. In this case, $\bs_{l} = \bs_{l-3}zyz = \bs_{l-1}z$, and so 
$\bs_{m}^{(1)} = \bs_{m-1}^{(1)}y^\prime = \bs_{m-2}^{(1)}yy^\prime$ for some $y^\prime \in \cA$.  If $z \ne a$, then $zw = za\bs_{l-3}zx$ is a factor of $\bt$ since $w \in F(\bt)$ and $z$ is 
separating for $\bt$. So, letting $w^\prime = zw$ if $z \ne a$ and $w^\prime = w$ if $z = a$, we have 
$w^\prime \in F(\bt)$ and
\begin{align*} 
\Psi_z^{-1}(w^\prime) &= a\Psi_{z}^{-1}(\bs_{l-3}zx)   \\
&= a\Psi_{z}^{-1}(\bs_{l-1}y^{-1}x) \\
&= a\bs_{m-1}^{(1)}(z^{-1}y)^{-1}\Psi_z^{-1}(x) \\
&=a\bs_{m-1}^{(1)}y^{-1}z\Psi_z^{-1}(x) \\
&=a\bs_{m-2}^{(1)}z\Psi_z^{-1}(x).
\end{align*}
That is, $w^\prime = \Psi_z(w^{(1)})$ where $w^{(1)} \in F(\bt^{(1)})$ is given by
\[
 w^{(1)} = \begin{cases} 
                  a\bs_{m-2}^{(1)}xx &\mbox{if $x = z$}, \\
                  a\bs_{m-2}^{(1)}x &\mbox{if $x \ne z$}.
                  \end{cases}
\]
Therefore, since $\bs_{m}^{(1)} = \bs_{m-1}^{(1)}y^\prime = \bs_{m-2}^{(1)}yy^\prime$ with $y> x$, we have  $w^{(1)} < a\bs_{m-1}^{(1)}$; a contradiction. 

Both Cases 1 and 2 lead to a contradiction; whence $\min(\bt) = a\bs$. 
\qed\end{pf*}\vspace{-10pt}

We now prove our main result: a characterization of fine words over a finite alphabet. (Recall that $\bv_p$ denotes the prefix of length $p$ of a given infinite word $\bv$, and $\rev{\bv}_p$ denotes its reversal.)

\begin{thm} \label{T:fine} Suppose $\bt$ is an infinite word with $\mbox{{\em Alph}}(\bt) = \cA$. Then, $\bt$ is fine if and only if one of the following holds:
\begin{enumerate}
\item[i)] $\bt$ is a strict episturmian word; 
\item[ii)] $\bt = v\mu(\bv)$ where $\bv$ is a $\cB$-strict standard episturmian word with $\cB = \cA\setminus\{x\}$, $\mu$ is a pure epistandard morphism on $\cA$,  and $v$ is a non-empty suffix of $\mu(\rev{\bv}_px)$ for some $p \in \NN$.  
\end{enumerate}
\end{thm}
\vspace{-15pt}\begin{pf}
In what follows, let $\cA$ denote the alphabet of $\bt$. 

\noindent ($\Rightarrow$): $\bt$ is fine, so there exists an infinite word $\bs$ such that, for any letter $a \in \cA$ and lexicographic order $<$ satisfying $a = \min(\cA)$, we have $\min(\bt) = a\bs$. Further, $\min(\bt) \leq \min(\bs)$ since $F(\bs) \subseteq F(a\bs) \subseteq F(\bt)$, and therefore $a\bs \leq \min(\bs)$.  
Thus, Proposition \ref{P:gP-2} implies that $\bs$ is a standard episturmian word over $\cA$.  
We distinguish two cases, below.

\noindent \emph{Case $1$}: $a\bs = \min(\bs)$ for any $a \in \cA$ and lexicographic order such that $a = \min(\cA)$. 

By Proposition \ref{P:jJgP-1}, $\bs$ is an $\cA$-strict standard episturmian word.

Clearly, $F(\bs) \subseteq F(\bt)$ and we show that $F(\bs) = F(\bt)$ which implies $\bt$ is equivalent to $\bs$, and hence $\bt$ is an $\cA$-strict episturmian word. Suppose, on the contrary, $F(\bs) \ne F(\bt)$. Then there exists a word $u \in F(\bt)\setminus F(\bs)$, say $u = xv$ with $|v| = m$ minimal and $x \in \cA$. Now, $v$ is not a prefix of $\bs$; otherwise, $zv$ is a factor of $\bs$ for all $z \in \cA$ (since any prefix of $\bs$ is left special and has $|\cA|$ distinct left extensions in $\bs$), which contradicts the fact that $xv \not\in F(\bs)$. Therefore, for some order such that $x = \min(\cA)$, we have $u=xv < x\bs_m$, contradicting the fact that $\min(\bt) = x\bs$; whence $F(\bt) = F(\bs)$. 

\noindent \emph{Case $2$}: $x\bs < \min(\bs)$ for some $x \in \cA$ and lexicographic order such that $x = \min(\cA)$. 

In this case, it follows from Propositions \ref{P:jJgP-1} and \ref{P:gP-2}  that $\bs \in \cAw$ is a standard episturmian word that is not $\cA$-strict. Therefore,  letting $\Delta(\bs) = x_1x_2x_3\cdots$, there exists a minimal $n \in \NN$ such that  
$\bs = \Psi_{x_1}\Psi_{x_2}\cdots\Psi_{x_n}(\bs^{(n)}) = \mu_n(\bs^{(n)})$, where $\bs^{(n)}$ is a $\cB$-strict standard episturmian word with $\cB =$~Alph$(\Delta(\bs^{(n)}))=$~Alph$(x_{n+1}x_{n+2}x_{n+3}\cdots) \subset \cA$ (see Remark \ref{R:ult_remark}).
Note that if $\bs = \bs^{(0)}$ is $\cB$-strict with $\cB \subset \cA$, then $n=0$ and we take $x_n = x_0$ to be a letter in $\cA\setminus \cB$.

Clearly, $\bs$ begins with $x_1 \in \cA$ and $x_1$ is separating for $\bs$. Observe that $x_1$ must also  be separating for $\bt$. Indeed, let us suppose that this is not true. Then, there exist letters $z$, $z^\prime \in \cA\setminus \{x_1\}$ (possibly equal) such that $zz^\prime \in F(\bt)$. But, if $<$ is an  order such that $\min(\cA) = z \leq z^\prime < x_1$, then $zx_1$ is a prefix of $z\bs$ with $zz^\prime <  zx_1$,  contradicting the fact that $\min(\bt) = z\bs$. Therefore, $x_1$ must be separating for $\bt$.

Now, let $<$ be an order with $a = \min(\cA)$. Let $\bt^\prime = \bt$ if $\bt$ begins with $x_1$. Otherwise, if $\bt$ begins with $y \ne x_1$, let $\bt^\prime = x_1\bt$. In the latter case, $ax_1 \pref a\bs$ and $x_1y \pref \bt^\prime$ with $ax_1 < x_1y$; thus $\min(\bt^\prime) = \min(\bt) = a\bs$.       
So we may consider $\bt^\prime$ instead of $\bt$. 

Observe that $\bs = \Psi_{x_1}(\bs^{(1)})$ and, since $x_1$ is separating for $\bs$ (and hence for $\bt^\prime$), we have $\bt^\prime = \Psi_{x_1}(\bt^{(1)})$ for some $\bt^{(1)} \in \cA^\omega$. Because $\min(\bt^\prime) = x\bs$ for any letter $x \in \cA$ and lexicographic order such that $x = \min(\cA)$, it follows that Alph$(\bt^{(1)}) = \cA$ (see  arguments in the first lines of the proof of Lemma \ref{L:25.01.2006}); in particular $x_1 \in \mbox{Alph}(\bt^{(1)})$. So, by  Lemma \ref{L:25.01.2006}, we have 
\[
\min(\bt^{(1)}) = a\bs^{(1)}. 
\]
Continuing in the same way (and applying Lemma \ref{L:25.01.2006} repeatedly), we obtain sequences $(\bs^{(i)})$, $(\bt^{(i)})$, $({\bt^\prime}^{(i)})$ for $i = 0, 1, 2$, $\ldots$~, $n$ such that $\bs^{(i-1)} = \Psi_{x_i}(\bs^{(i)})$, ${\bt^\prime}^{(i-1)} = \Psi_{x_i}(\bt^{(i)})$, where ${\bt^\prime}^{(i-1)} = \bt^{(i-1)}$ if $\bt^{(i-1)}$ begins with $x_i$, ${\bt^\prime}^{(i-1)} = x_i\bt^{(i-1)}$ otherwise, and $\bt^{(0)} = \bt$, ${\bt^\prime}^{(0)} = \bt^\prime$. In particular, we have Alph$(\bt^{(n)}) = \cA$ and  
\begin{equation*}
 \min(\bt^{(n)})  = a\bs^{(n)} 
\end{equation*}
for any $a \in \cA$ and lexicographic order $<$ satisfying $a = \min(\cA)$.
   
Now we show that $\cB = \cA\setminus\{x_n\}$, i.e., $\cA = \cB\cup\{x_n\}$. First observe that $x_n \in \cA\setminus\cB$ by minimality of $n$. 

Suppose $\bt^{(n)}$ contains two occurrences of the letter $x_n$. Then, since $x_{n+1}$ is separating for $\bt^{(n)}$, we have $x_nw^{(n)}x_n \in F(\bt^{(n)})$ for some non-empty word $w^{(n)}$ for which $x_{n+1}$ is separating, and the first and last letter of $w^{(n)}$ is $x_{n+1}$ (that is, $w^{(n)}x_n = \Psi_{x_{n+1}}(w^{(n+1)}x_n)$, where $w^{(n+1)} = \Psi_{x_{n+1}}^{-1}(w^{(n)}x_{n+1}^{-1})
$). Continuing the above procedure, we obtain infinite words $\bt^{(n+1)}$, $\bt^{(n+2)}$, $\ldots$ containing similar shorter factors $x_nw^{(n+1)}x_n$, $x_nw^{(n+2)}x_n$, $\ldots$ until we reach $\bt^{(q)}$, which contains $x_nx_n$. But this is impossible because $x_{q+1} \in \cB \subseteq \cA\setminus\{x_n\}$ is separating for $\bt^{(q)}$. Therefore, $\bt^{(n)}$ contains only one occurrence of $x_n$ and we have 
$$\bt^{(n)} = ux_n\bv \quad \mbox{for some $u \in (\cA\setminus\{x_n\})^*$ and $\bv\in (\cA\setminus\{x_n\})^\omega$}.$$ 
Note that the same reasoning allows to prove the unicity of $x_0 \in \cA\setminus\cB$ when $n=0$.

Clearly, for any order such that $x_n = \min(\cA)$, we have $\min(\bt^{(n)}) = x_n\bv = x_n\bs^{(n)}$; whence $\bv = \bs^{(n)}$ and so $\bt^{(n)} = ux_n\bs^{(n)}$. Note that if $u\ne \empt$, then $u$ ends with $x_{n+1}$, and in particular $x_{n+1}$ is separating for $ux_n$ since $x_{n+1}$ is separating for $\bt^{(n)}$. 
 
Let $u^\prime = x_{n+1}u$ if $u$ does not begin with $x_{n+1}$; otherwise let $u^\prime = u$. Then $u^\prime x_n$ is a prefix of ${\bt^{\prime}}^{(n)}$. Moreover, since $x_{n+1}$ is separating for $u^\prime x_{n}$,  we have $u^\prime x_{n} = \Psi_{x_{n+1}}(u^{(n+1)}x_n)$ where 
$u^{(n+1)} = \Psi_{x_{n+1}}^{-1}(u^\prime x_{n+1}^{-1})$. Hence $\bt^{(n+1)} = u^{(n+1)}x_n\bs^{(n+1)}$, where $x_{n+2}$ is separating for $u^{(n+1)}x_n$ (if $u^{(n+1)} \ne \empt$).  Continuing in this way,  we arrive at the infinite word $\bt^{(q)} = x_n\bs^{(q)}$ for some $q \geq n$. 

Now, reversing the procedure, we find that 
\[
  \bt^{(n)}  
  = w\bs^{(n)} \quad \mbox{where $w = ux_n$ is a non-empty suffix of $\Psi_{x_{n+1}}\cdots \Psi_{x_q}(x_n)$}.
\]
Accordingly, $u \in \cB^*$ since $x_{n+1}$, $\ldots$~, $x_{q} \in \cB$; whence $\cA = \cB \cup \{x_n\}$. 

Suppose $(u_i)_{i\geq1}$ is the sequence of palindromic prefixes of $\bs$ and the words $(h_i)_{i\geq0}$ are the prefixes  $(\mu_{i}(x_{i+1}))_{i\geq 0}$ of $\bs$. Then, letting $u_i^{(n)}$,  $h_i^{(n)}$, and $\mu_i^{(n)}$ denote the analogous elements for $\bs^{(n)}$, we have  
$$\mu_0^{(n)} = \mbox{Id}, \quad \mu^{(n)}_{i} = \Psi_{x_{n+1}}\Psi_{x_{n+2}}\cdots \Psi_{x_{n+i}} = \mu_n^{-1}\mu_{n+i}$$ and 
 \[
  h_0^{(n)} = x_{n+1}, \quad h_i^{(n)} = \mu_{i}^{(n)}(x_{n+1+i}) \quad \mbox{for $i=1$, $2$, $\ldots$~.}
\]
Now, if $u\ne \empt$, then $q \geq n+1$, and we have 
\begin{align*}
 \Psi_{x_{n+1}}\cdots \Psi_{x_q}(x_n) = \mu_{q-n}^{(n)}(x_n)  
                        &= \mu_{q-n-1}^{(n)}\Psi_q(x_n) \notag\\
                        &= \mu_{q-n-1}^{(n)}(x_qx_n) \notag\\
                        &= h_{q-n-1}^{(n)}\mu_{q-n-1}^{(n)}(x_n) \notag\\
                        &\qquad \vdots \notag\\
                        &=h_{q-n-1}^{(n)}\cdots h_{1}^{(n)}\mu_0^{(n)}(x_{n+1}x_n) \notag \\
                        &=h_{q-n-1}^{(n)}\cdots h_1^{(n)}h_0^{(n)}x_n  
                        =u_{q-n+1}^{(n)}x_n \quad \mbox{(by \eqref{eq:u_n}).} 
\end{align*}  
Therefore, $w = ux_n$ where $u$ is a (possibly empty) suffix of the palindromic prefix $u_{q-n+1}^{(n)}$  of $\bs^{(n)}$.  That is, $u$ is the reversal of some prefix of $\bs^{(n)} = \bv$; in particular 
$$u = \rev{\bv}_p \quad \mbox{for some $p \in \NN$},$$
and hence
\[
  \bt^{(n)} = \rev{\bv}_px_n\bv. 
\]  
So, passing back from $\bt^{(n)}$ to $\bt$, we find that   
\[
  \bt = v\mu_n(\bv) = v\bs \quad \mbox{where $v$ is a non-empty suffix of $\mu_n(\rev{\bv}_px_n)$.}
\]

Cases 1 and 2 give properties $i)$ and $ii)$, respectively.
\bigskip
   
\noindent ($\Leftarrow$):  Firstly, if $\bt$ is an $\cA$-strict episturmian word, then Proposition \ref{P:jJgP-1} implies that $\bt$ is fine. 

Now suppose $\bt = v\mu(\bv)$ where  $\bv$ is a $\cB$-strict standard episturmian word with $\cB = \cA\setminus\{x\}$, $\mu$ is a pure epistandard morphism on $\cA$, and $v$ is a non-empty suffix of $\mu(\rev{\bv}_px)$ for some $p \in \NN$. First observe that if $\mu =$ Id, then  
\[
 \bt = \rev{\bv}_qx\bv \quad \mbox{for some $q \leq p$.}
\] 
 
Consider an order $<$ such that $\min(\cA) = a \ne x$. Then, by Proposition \ref{P:jJgP-1}, $\min(\bv) = a\bv$, and it follows that $\min(\bt) = \min(\bv) = a\bv$. Indeed, if $\bt = x\bv$ (i.e., $q=0$), then it is clear that $\min(\bt) = \min(\bv)$. On the other hand, if $q\geq1$, let us suppose, on the contrary, that $\min(\bt) \ne \min(\bv)$. Then $\min(\bt)$ is a suffix of $\rev{\bv}_qx\bv$ containing the letter $x$, i.e., 
\[
  \min(\bt) = a\rev{\bv}_lx\bv \quad \mbox{for some $l$ with $0 \leq l < q \leq p$}
\]
(where $a\rev{\bv}_l = \rev{\bv}_{l+1} \in F(\bv)$). But, since $\min(\bv) = a\bv$, $a\bv_{l+1}$ is a factor of $\bv$ (and hence a factor of $\bt$) with $a\bv_{l+1} = a\bv_l a < a\rev{\bv}_lx$; a contradiction. Thus  $\min(\bt) = \min(\bv) = a\bv$. Moreover, it is clear that $\min(\bt) = x\bv$ for any order such that  $x = \min(\cA)$.  So we have shown that  
$\min(\bt) = a\bv$ for any letter $a\in \cA$ and lexicographic order satisfying $a = \min(\cA)$; whence $\bt$ is fine.   

Now consider the case when $\mu$ is not the identity. Let us suppose that $\bt$ is not fine and let $\mu$ be minimal with this property. Then, $\mu = \Psi_z\eta$ for some $z \in \cA$ and pure epistandard morphism $\eta$.  

Consider $\bt^\prime = v^\prime\mu(\bv)$, where $v^\prime = v$ if $v$ begins with $z$ and $v^\prime = zv$ otherwise. Then $v^\prime$ is also a non-empty suffix of $\mu(\rev{\bv}_px)$ since $z$ is separating for the word $\mu(\rev{\bv}_px)$ (which begins with $z$). Letting $\bt^{(1)} = \Psi_z^{-1}(\bt^\prime)$, we have   
\[
  \bt^{(1)} = w\eta(\bv) 
\]
where $w = \Psi_z^{-1}(v^\prime)$ is a non-empty suffix of $\eta(\rev{\bv}_px)$.    
By minimality of $\mu$, $\bt^{(1)}$ is fine, so there exists an infinite word $\bs^{(1)} = \Psi_z^{-1}(\bs)$ such that $\min(\bt^{(1)}) = a\bs^{(1)}$ for any $a \in \cA$ and order $<$ satisfying $a = \min(\cA)$. But then $\min(\bt) = \min(\bt^\prime) = a\bs$ by Lemma \ref{L:25.01.2006}. Thus $\bt$ is fine.
\qed\end{pf}\vspace{-10pt}

\begin{example} Let $\cA = \{a,b,c\}$ with $a<b<c$ and suppose $\bbf$ is the Fibonacci word over $\{a,b\}$. Then, the following infinite words are fine. 
\begin{itemize}
\item $\bbf = abaababaabaaba\cdots$ \smallskip
\item $c\bbf = \underline{c}abaababaabaaba\cdots$ \smallskip
\item $\rev{\bbf_4}c\bbf = aaba\underline{c}abaababaabaaba\cdots$ \smallskip
\item $\Psi_a(\bbf) = aabaaabaabaaabaaaba\cdots$ \smallskip
\item $\Psi_c(c\bbf) = \underline{c}cacbcacacbcacbcacacbcacacbca\cdots$ \smallskip 
\item $\Psi_c(\rev{\bbf_4}c\bbf) = cacacbca\underline{c}cacbcacacbcacbcacacbcaca\cdots$ \smallskip
\end{itemize}
Let us note, for example, that $\Psi_c(\bbf)$ is {\bf not} fine since it is a {\em non-strict} standard episturmian word. That is, $\Psi_c(\bbf)$ is a standard episturmian word with directive word $c(ab)^\omega$, so it is not strict, nor does it take the second form given in Theorem \ref{T:fine}. 
\end{example}

\section{Concluding remarks}

It is easy to see that Proposition \ref{P:gP05mors} is a special case of Theorem \ref{T:fine} because the 2-strict episturmian words are precisely the Sturmian words and the 1-strict standard episturmian words are periodic infinite words of the form $x^\omega$ where $x$ is a letter (see \cite[Proposition 2.9]{jJgP02epis}).  

As alluded to in the introduction, an infinite word taking form $ii)$ in Theorem \ref{T:fine} is said to be a strict {\em skew episturmian} word. Skew episturmian words (now called {\em episkew words} \cite{jAaG07extr, aGjJgP06char}) are explicated in the paper \cite{aGjJgP06char}, in which we expand on our work here by characterizing via lexicographic order all {\em episturmian words in a wide sense}, i.e., all infinite words whose factors are (finite) episturmian.

\ack{The author would like to thank Jacques Justin for suggesting the problem and giving many helpful comments on preliminary versions of this paper. Thanks also to the referees for their careful reading of the paper and providing thoughtful suggestions.}

\end{document}